\documentstyle[12pt,amscd]{amsart}
\newtheorem{Theorem}{Theorem}[section] \newtheorem{Lemma}[Theorem]{Lemma}
\newtheorem{Corollary}[Theorem]{Corollary}

\def\reg{\operatorname{reg}} \def\In{\operatorname{in}}  \def\gin{\operatorname{Gin}} \def\char{\operatorname{char}}   \def\To{\longrightarrow} \def\syz{\operatorname{Syz}}  \def\ext{\operatorname{Ext}}
\def\Hom{\operatorname{Hom}}

\begin{document}
\title{Gr\"obner bases, local cohomology and \\ reduction number}
\author{Ng\^o Vi\^et Trung}
\address{Institute of Mathematics,  Box 631, B\`o H\^o, Hanoi, Vietnam}
\email{nvtrung@@hn.vnn.vn}
\thanks{The author is partially supported by the National Basic Research Program}
\begin{abstract} D. Bayer and M. Stillman showed that Gr\"obner bases can be used to compute the Castelnuovo-Mumford regularity which is a measure for the vanishing of graded local cohomology modules. The aim of this paper is to show that the same method can be applied to study other cohomological invariants as well as the reduction number. \end{abstract}
\maketitle

\section*{Introduction} \smallskip
 
Let $S = k[x_1,\ldots,x_n]$ be a polynomial ring over a field $k$ of arbitrary characteristic. Let $\frak M$ be the maximal graded ideal of $S$. For any finitely generated graded $S$-module $M$  we will denote by $H_{\frak M}^i(M)$ the $i$th local cohomology of $M$ with respect to $\frak M$. Since $H_{\frak M}^i(M)$ is an artinian graded module, we may consider the largest non-vanishing degree  
$$a_i(M) = \max\{n|\, H_{\frak M}^i(M)_n \neq 0\}$$
with the convention $a_i(M) = -\infty$ if $H_{\frak M}^i(M) = 0$. Note that $H_{\frak M}^i(M) = 0$ for $i > \dim M$. The  Castelnuovo-Mumford regularity of $M$ is defined as
$$\reg(M) = \max\{a_i(M)+i|\ i \ge 0\}.$$
This invariant  carries important information on the structure of $M$ [EG], [O]. Similarly, we define
$$a^*(M) = \max\{a_i(M)|\ i \ge 0\}.$$
It is known that $a^*(M)+1$ gives an upper bound for the regularity of the Hilbert function of $M$.  Moreover, $a^*(M)$ can be used to estimate $-a_d(M)$, $d = \dim M$, which is equal to the least non-vanishing degree of the canonical module of $M$ [GW]. The Castelnuovo-Mumford regularity $\reg(M)$ and the largest non-vanishing degree $a^*(M)$ of local cohomology modules can be viewed as special cases of the more general invariants:
\begin{eqnarray*} \reg_t(M) &  = & \max\{a_i(M)+i|\ i \le t\},\\
a^*_t(M) & = & \max\{a_i(M)|\ i \le t\},\end{eqnarray*}
where $t = 0,\ldots,d$. These invariants have been studied in [T1], [T2], [T3].\par

For any homogeneous ideal $I \neq 0$ in $S$ we have $\reg_t(S/I) = \reg_t(I)-1$ and $a^*_t(S/I) = a^*_t(I)$. Let $\gin(I)$ denote the generic initial ideal of $I$  with respect to a given term order  in generic coordinates. Bayer and Stillman [BS2] proved that $\reg(I) = \reg(\gin(I))$ with respect to the reverse lexicographic term order and that if $\char(k) = 0$, then $\reg(\gin(I))$ is the maximum degree of the minimal generators of $\gin(I)$. We will use their method to prove the following similar statements on $\reg_t(I)$ and  $a^*_t(I)$. \medskip

\noindent {\bf Corollary 1.4. } {\em  Let $\gin(I)$ denote the generic initial ideal of $I$ with respect to the reverse lexicographic order. Then\par \smallskip
{\rm (i)} $\reg_t(I) = \reg_t(\gin(I))$,\par \smallskip
{\rm (ii)} $a^*_t(I) = a^*_t(\gin(I))$.}\medskip

\noindent{\bf Corollary 2.5. }  {\em Assume that $\char(k) = 0$. For any monomial $x^A$ let $m(x^A)$ denote the largest $i$ such that $x_i$ divides $x^A$. Then\par \smallskip
{\rm (i)} $\reg_t(\gin(I))$ is the maximum degree of the minimal generators $x^A$ of $\gin(I)$ with $m(x^A) \ge n-t$,\par \smallskip
{\rm (ii)} $a_t^*(\gin(I))$ is the maximum of $\deg(x^A)+m(x^A)-n$ of the minimal generators $x^A$ of $\gin(I)$ with $m(x^A) \ge n-t$.}\medskip

The equality $\reg_t(I) = \reg_t(\gin(I))$ was already  proved by D. Bayer, H. Charalambous, and S. Popescu [BCP] from a different point of view.  Let
$$0 \To F_s \To \ldots \To F_1 \To F_0 \To M$$
be a minimal free resolution of $M$ over $S$.  Write $b_i$ for the maximum degree of the generators of $F_i$. Then the Castelnuovo-Mumford regularity $\reg(M)$ can be also defined by the formula $\reg(M) = \max\{b_i-i| i \ge 0\}.$ Motivated by this definition Bayer, Charalambous, and Popescu introduced the $l$-regularity of $M$ as $$l\text{-}\reg(M) = \max\{b_i-i| i \ge l\}$$
and showed that $l$-$\reg(I) = l$-$\reg(\gin(I))$ for the reverse lexicographic order. We shall see that
\begin{eqnarray*} \reg_t(M) & = & \max\{b_i-i|\ i \ge n-t\},\\
a^*_t(M) & = & \max\{b_i|\ i \ge n-t\}. \end{eqnarray*}
Therefore, $\reg_t(M) = (n-t)$-$\reg(M)$. It  should be pointed out that Bayer, Charalambous and Popescu proved more  than the equality $l$-$\reg(I) = l$-$\reg(\gin(I))$, namely that the extremal Betti numbers of $I$  which correspond to the  ``jumps" in the regularity of the successive syzygy modules do not change when passing to a generic initial ideal of $I$. This result was extended to exterior algebras by A. Aramova and J. Herzog [AH2].  \par

We will also use the method of Bayer and Stillman to study the reduction number of  a graded algebra. Let $\frak m$ denote the maximal graded ideal of  $S/I$.  An ideal $\frak a$ of $S/I$ is called a reduction of $\frak m$ if $\frak m^{r+1} = \frak a\frak m^r$ for large $r$. The least number $r$ with this property is denoted by $r_{\frak a}(S/I)$. A reduction of $\frak m$ is said to be minimal if it does not contain any other reduction of $\frak m$. The reduction number $r(S/I)$ of $S/I$ is defined as the minimum $r_{\frak a}(S/I)$ of all minimal reductions $\frak a$ of $\frak m$. This number may be used to estimate $\reg(S/I)$ and $a^*(S/I)$ [T1], [T2], [T3]. One may view $r(S/I)$ as a measure for the complexity of $S/I$ [V]. However, the relationship between $r(S/I)$ and $r(S/\In(I))$ is not well-understood. W. Vasconcelos conjectured that $r(S/I) \le r(S/\In(I))$ [V, Conjecture 5.15]. Recently, Bresinsky and Hoa [BH] proved that $r(S/I) \le r(S/\gin(I))$. Inspired by their result  we will show that equality holds for the reverse lexicographic order. \medskip

\noindent{\bf Theorem 4.3. } {\em Let $\gin(I)$ denote the generic initial ideal of $I$ with respect to  the reverse lexicographic order. Assume that $k$ is an infinite field. Then $$r(S/I) = r(S/\gin(I)).$$} 

This theorem gives a pratical way to compute $r(S/I)$ by means of Gr\"obner bases since $r(S/\gin(I))$  is the least number $r$ for which $x_{n-d}^{r+1} \in \gin(I)$ [BH]. We also give an example showing that Theorem 4.3 does not hold for an arbitrary term order. \par

The paper is divided into four sections. Section 1 deals with the invariance of $\reg_t(I)$ and $a^*_t(I)$ when passing to certain initial ideal of $I$. Section 2 gives combinatorial descriptions of $\reg_t(I)$ and $a^*_t(I)$ when $I$ is a Borel-fixed ideal. Section 3 investigates the relationships between the invariants $\reg_t(M)$ and $a^*_t(M)$ and the syzygies of a graded module $M$. Section 4 is devoted to the study of the reduction number of $S/I$. For unexplained terminology we refer to the book of Eisenbud [E].\smallskip

\noindent{\em Acknowledgement. }  The author is grateful to J. Herzog  and L. T. Hoa for helpful suggestions. 

\section{Gr\"obner bases and cohomological invariants} \smallskip

Let $I$ be an arbitrary homogeneous ideal in the polynomial ring $S = k[X]$. Let $R$ denote the factor ring $S/I$. It is known that the cohomological invariants $\reg_t(R)$ and $a^*_t(R)$ can be characterized in terms of a sequence $z_1,\ldots,z_{t+1}$ of linear forms in $R$.\par

Recall that  $z_1,\ldots,z_s$ of homogeneous elements in $R$ is called a {\it filter-regular sequence} in $R$ if $z_i \not\in \frak p$ for any associated prime ideal $\frak p \neq \frak m$ of $(z_1,\ldots,z_{i-1})$, $i = 1,\ldots,s$, where $\frak m$ denotes the maximal graded ideal of $R$. Note that if $(z_1,\ldots,z_{i-1})$ has no associated prime ideal $\frak p \neq \frak m$, then $(z_1,\ldots,z_{i-1})$ is an $\frak m$-primary ideal and $z_i$ can be any homogeneous element of $R$. If $k$ is an infinite field, we may always assume that $x_1,\ldots,x_n$ is a filter-regular sequence in $R$ by a generic choice of variables. \par

Filter-regular sequences have their origin in the theory of Buchsbaum rings [STC].  Here we are mainly interested in the criterion that $z_1,\ldots,z_s$ is a filter-regular sequence if and only if
$$[(z_1,\ldots,z_{i-1}):z_i ]_m = (z_1,\ldots,z_{i-1})_m,\; i = 1,\ldots,s,$$
for large $m$  [T1, Lemma 2.1]. Such a sequence of linear forms was named almost regular in [AH2].

\begin{Theorem}  {\rm [BS2, Theorem (1.10)], [T1, Proposition 2.2], [T3, Corollary 2.6]} Let $z_1,\ldots,z_{t+1}$ be  a filter-regular sequence of linear forms in $R$. Then\par \smallskip
{\rm (i) } $\reg_t(R)$ is the largest integer $r$ for which there is an index $0 \le i \le t$ such that
$$[(z_1,\ldots,z_i):z_{i+1}]_r \neq (z_1,\ldots,z_i)_r,$$ \par  
{\rm (ii) } $a_t^*(R)$ is the largest interger $a$ for which there is an index $0\le i \le t$ such that
$$[(z_1,\ldots,z_i):z_{i+1}]_{a+i} \neq (z_1,\ldots,z_i)_{a+i}.$$\end{Theorem}

The above characterizations of $\reg_t(R)$ and $a^*_t(R)$ provide a link to Gr\"obner bases by means of  the following result of Bayer and Stillman. 

\begin{Lemma}  {\rm [BS2, Lemma (2.2)]}  Let $\In(I)$ denote the initial ideal with respect to  the reverse lexicographic order of $I$. Let $i = n,\ldots,1$. For every integer $m \ge 0$, $$[(I,x_n,\ldots,x_{i+1}):x_i]_m = (I,x_n,\ldots,x_{i+1})_m$$ if and only if  $$[(\In(I),x_n,\ldots,x_{i+1}):x_i]_m = (\In(I),x_n,\ldots,x_{i+1})_m.$$ \end{Lemma}

\begin{Theorem}   Let $I$ be an arbitrary homogeneous ideal. Let $\In(I)$ denote the initial ideal of $I$ with respect to the reverse lexicographic order. Assume that $x_n,\ldots,x_1$ is a filter-regular sequence in $S/I$. Then\smallskip \par
{\rm (i) } $\reg_t(I) = \reg_t(\In(I))$,\smallskip \par
{\rm (ii) } $a_t^*(I) = a_t^*(\In(I))$.\end{Theorem}

\begin{pf}   For any graded $S$-module $M$ let $\delta(M)$ denote the largest integer $r$ such that $M_r \neq 0$ with the convention $\delta(M) = -\infty$ if $M = 0$ and $\delta(M) = \infty$ if $M$ is not of finite length. For $i = n,\ldots,1$, $\delta\big((I,x_n,\ldots,x_{i+1}):x_i/(I,x_n,\ldots,x_{i+1})\big)$ is just the largest integer $r$ such that
$[(I,x_n,\ldots,x_{i+1}):x_i]_r \neq (I,x_n,\ldots,x_{i+1})_r.$ The assumption that $x_n,\ldots,x_1$ is a filter-regular sequence in $S/I$ implies that $\delta\big((I,x_n,\ldots,x_{i+1}):x_i/(I,x_n,\ldots,x_{i+1})\big) < \infty$
for $i = n,\ldots,1$. By Lemma 1.2, 
\begin{eqnarray*} & \delta\big ((I,x_n,\ldots,x_{i+1}):x_i/(I,x_n,\ldots,x_{i+1})\big) \\ &  = \delta\big((\In(I),x_n,\ldots,x_{i+1}):x_i/(\In(I),x_n,\ldots,x_{i+1})\big).\end{eqnarray*}
Hence $\delta\big((\In(I),x_n,\ldots,x_{i+1}):x_i/(\In(I),x_n,\ldots,x_{i+1})\big) < \infty$ for $i = n,\ldots,1$. So $x_n,\ldots,x_1$ is a filter-regular sequence in $S/\In(I)$. Note that $\reg_t(I) = \reg_t(R)+1$ and $a^*_t(I) = a^*_t(R)$. Applying Theorem 1.1 we obtain
\begin{eqnarray*} \reg_t(I)& =& \max\big\{\delta((I,x_n,\ldots,x_{i+1}):x_i/(I,x_n,\ldots,x_{i+1}))|\, i =n,\ldots,n-t\big\}+1\\
&=&\max\{\delta\big((\In(I),x_n,\ldots,x_{i+1}):x_i/(\In(I),x_n,\ldots,x_{i+1})\big)|\, i = n,\ldots,n-t\big\}+1\\
&= & \reg_t(\In(I)). \end{eqnarray*}
Similarly we have
\begin{eqnarray*} a^*_t(I)& =& \max\big\{\delta\big((I,x_n,\ldots,x_{i+1}):x_i/(I,x_n,\ldots,x_{i+1})\big)-n+i|\, i =n,\ldots,n-t\big\}\\
&= &\max\{\delta\big((\In(I),x_n,\ldots,x_{i+1}):x_i/(\In(I),x_n,\ldots,x_{i+1})\big)-n+i|\, i = n,\ldots,n-t\big\}\\ &= & a^*_t(\In(I)).  \end{eqnarray*} \end{pf} 

Let the general linear group GL$(n,k)$ of invertible $n\times n$ matrices over $k$ act as a group of algebra automorphisms on $S = k[x_1,\ldots,x_n]$.  There is a Zariski open set $U \subset$ GL$(n,k)$ and a monomial ideal $J \subset S$ such that for all $g \in U$ we have $\In(gI) = J$ (see e.g. [E, Theorem 15.18]). The ideal $J$ is called a {\it generic initial ideal} of $I$, denoted by $\gin(I)$.

\begin{Corollary}  Let $\gin(I)$ denote the generic initial ideal with respect to the reverse lexicographic order. Then\smallskip \par
{\rm (i)} $\reg_t(I) = \reg_t(\gin(I))$,\smallskip\par
{\rm (ii)} $a^*_t(I) = a^*_t(\gin(I))$. \end{Corollary}

 \begin{pf}   For a generic choice of coordinates we may assume that $x_n,\ldots,x_1$ is a filter-regular sequence in $S/I$. Hence the conclusions follow from Theorem 1.3.  \end{pf}  

As mentioned in the introduction, the first statement of Corollary 1.4 can be deduced from a recent result of Bayer, Charalambous, and Popescu [BCP]. This will be discussed in Section 3. 
 
\section{Cohomological invariants of Borel-fixed monomial ideals}\smallskip

Let $\cal B$ be the Borel subgroup of GL$(n,k)$ consisting of the upper triangular invertible matrices. A monomial ideal $I$ is called {\it Borel-fixed} if for all $g \in \cal B$, $g(I) = I$.  

\begin{Theorem} {\rm [Ga], [BS1]}  Let $I$ be an arbitrary homogeneous ideal. Then $\gin(I)$ is a Borel-fixed ideal. \end{Theorem}

Borel-fixed ideals can be characterized as follows.

\begin{Lemma} {\rm  [BS2, Proposition 2.7]}  Let $I$ be a monomial ideal. Assume that $\char(k) = 0$. Then $I$ is Borel-fixed if and only if whenever $x_1^{p_1}\cdots x_n^{p_n} \in I$, then 
$$x_1^{p_1}\cdots x_i^{p_i+q}\cdots x_j^{p_j-q}\cdots x_n^{p_n} \in I$$
for each $1 \le i < j \le n$ and $0 \le q \le p_j$. \end{Lemma}

In the following we will denote a monomial of $S$ by $x^A$ and by $m(x^A)$ the largest $i$ such that $x_i$  divides $x^A$. 

\begin{Lemma}  Let $I$ be a Borel-fixed monomial ideal. Assume that $\char(k) = 0$. For $i = 1,\ldots,n$, let $r_i$ denote the largest integer $r$ such that 
$$[(I,x_n,\ldots,x_{i+1}):x_i]_r \neq (I,x_n,\ldots,x_{i+1})_r.$$
Then $r_i = \max\{\deg(x^A)|\ x^A\; \text{is a minimal generator of $I$ with $m(x^A)= i$}\}-1$. \end{Lemma}

 \begin{pf}  Let $r  = \max\{\deg(x^A)|\  \text{$x^A$ is a minimal generator of $I$ with $m(x^A) = i$}\}-1$. We fix a minimal generator  $x^A$ of $I$ of degree $r+1$ with $m(x^A) = i$. Write $x^A = x^Bx_i$. Then $x^B \in [(I,x_n,\ldots,x_{i+1}):x_i]_r$. Since $x^B \not\in I$ and since $x^B$ is not divisible by the variables $x_n,\ldots,x_{i+1}$, $x^B \not\in (I,x_n,\ldots,x_{i+1})$. Hence 
$$[(I,x_n,\ldots,x_{i+1}):x_i]_r \neq (I,x_n,\ldots,x_{i+1})_r.$$
To show that $r_i = r$ it suffices to show that
$$[(I,x_n,\ldots,x_{i+1}):x_i]_m = (I,x_n,\ldots,x_{i+1})_m,$$
for $m \ge r+1$. Assume to the contrary that there is a monomial $x^C \in [(I,x_n,\ldots,x_{i+1}):x_i]_m$ but $x^C \not\in (I,x_n,\ldots,x_{i+1})$.  Then $x^Cx_i \in (I,x_n,\ldots,x_{i+1})$ and $x^C$ is not divisible by the variables $x_n,\ldots,x_{i+1}$. Hence $x^Cx_i \in I$ and $m(x^Cx_i) =  i$. Since $\deg x^Cx_i = m+1 \ge r+2$, $x^Cx_i$ is not a minimal generator of $I$. Therefore we can find a monomial $x^D \in I$ such that $x^Cx_i = x^Dx_h$ for some $h \le i$. Since $x^C \not\in I$, $x^C \neq x^D$ so that $h \neq i$. Thus, $x^D$ is divisible by $x_i$ and we may write $x^D = x^Ex_i$. It follows that $x^C = x^Ex_h$. By Lemma 2.2, this implies $x^C \in I$, a contradiction. \end{pf} 

Now  we can describe the invariants $\reg_t(I)$ and $a_t^*(I)$ of a Borel-fixed ideal $I$ in terms of the minimal generators of $I$.  

\begin{Theorem}   Let $I$ be a Borel fixed monomial ideal.  Assume that $\char(k) = 0$. For any monomial $x^A$ in $S$ we denote by $m(x^A)$ the maximum of the index $j$ such that $x^A$ is divided by $x_j$. Then\smallskip\par
 {\rm (i)} $\reg_t(I)$ is the maximum degree of the minimal generators $x^A$ of $I$ with $m(x^A) \ge  n-t$,\smallskip\par
{\rm (ii)} $a_t^*(I)$ is the maximum of $\deg(x^A)+m(x^A)-n$ of the minimal generators $x^A$ of $I$ with $m(x^A) \ge n-t$. \end{Theorem}

\begin{pf}   For $i = n,\ldots,1$, let $r_i$ denote the largest integer $r$ such that 
$$[(I,x_n,\ldots,x_{i+1}):x_i]_r \neq (I,x_1,\ldots,x_{i+1})_r.$$
Then $r_i < \infty$ by Lemma 2.3. Hence $x_n,\ldots,x_1$ is a filter-regular sequence in $S/I$. Note that $\reg_t(I) = \reg_t(R)+1$ and $a^*_t(I) = a^*_t(R)$. By Theorem 1.1 and Lemma 2.3 we obtain
\begin{eqnarray*} \reg_t(I) &  = & \max\{r_i+1|\ i  = n,\ldots,n-t\}\\
& = & \max\{\deg(x^A)|\ x^A\; \text{is a minimal generator of $I$ with $m(x^A) \ge n-t$}\},\\
a^*_t(I) & = & \max\{r_i-n+i+1|\ i = n,\ldots,n-t\}\\
& = & \max\{\deg(x^A)+m(x^A)-n|\ x^A\; \text{is a minimal generator of $I$} \\ &&\quad\quad\quad \text{ with $m(x^A) \ge n-t$}\}.  \end{eqnarray*}\end{pf}

\noindent{\em Remark. } J. Herzog has informed the author that Theorem 2.4 can be derived from Eliahou-Kervaire's resolution for a stable monomial ideal [EK] (see also [AH1]). 

\begin{Corollary} Let $I$ be an arbitrary homogeneous ideal. Let $\gin(I)$ denote the generic initial ideal of $I$. Assume that $\char(k) = 0$. Then\smallskip\par
{\rm (i)} $\reg_t(\gin(I))$ is the maximum degree of the minimal generators $x^A$ of $\gin(I)$ with $m(x^A) \ge n-t$,\smallskip\par
{\rm (ii)} $a_t^*(\gin(I))$ is the maximum of $\deg(x^A)+m(x^A)-n$ of the minimal generators $x^A$ of $\gin(I)$ with $m(x^A) \ge n-t$. \end{Corollary}

\begin{pf}  By Theorem 2.1, $\gin(I)$ is a Borel-fixed ideal. Hence the conclusions follow from Theorem 2.4.  \end{pf} 

\section{Syzygies and cohomological invariants}\smallskip

Let $M$ be an arbitrary  graded module over the polynomial ring $S = k[x_1,\ldots,x_n]$. Let
$$0 \To F_s \To \ldots \To F_1 \To F_0 \To M$$
be a minimal free resolution of $M$ over $S$.  Write $b_i$ for the maximum degree of the generators of $F_i$. Motivated by the well-known formula $\reg(M) = \max\{b_i-i|\ i \ge 0\}$ Bayer, Charalambous and Popescu  [BCP] introduced the $l$-regularity 
$$l\text{-}\reg(M) = \max\{b_i-i|\ i \ge l\}$$ 
and proved that $l$-$\reg(I) =$ $l$-$\reg(\gin(I))$ for the reverse lexicographic order. Following an argument of Eisenbud in [E] we obtain the following relationships between the degree $b_i$ and the invariants $\reg_t(M)$ and $a^*_t(M)$. From this one can see that  
$\reg_t(M) = (n-t)\text{-}\reg(M).$  Hence the equality $\reg_t(I) = \reg_t(\gin(I))$ is only a consequence of the result of Bayer, Charalambous, and Popescu.

\begin{Theorem} Let $M$ be an arbitrary graded $S$-module of finite type. Then\smallskip\par
{\rm (i)}  $\reg_t(M) = \max\{b_i-i|\ i \ge n-t\}$,\smallskip\par
{\rm (ii)} $a^*_t(M) =  \max\{b_i|\ i \ge n-t\} - n$. \end{Theorem}

 \begin{pf}  By local duality (see e.g. [E, Theorem A4.2])  we have $H_{\frak M}^i(M)  = \ext_S^{n-i}(M,S(n))^{\vee},$ where $^\vee$ denotes the Matlis duality. From this it follows that 
$$a_i(M) = \max\{m|\ \ext_S^{n-i}(M,S)_{-m-n} \neq 0\}.$$ Hence
\begin{eqnarray*} \reg_t(M) & = & \max\{m|\ \ext_S^{n-i}(M,S)_{-m-n+i} \neq 0\; \text{for some}\; i \le t\}\\
& = & \max\{m|\  \ext_S^i(M,S)_{-m-i} \neq 0\; \text{for some}\; i \ge n-t\},\\
a^*_t(M) & = & \max\{m|\ \ext_S^{n-i}(M,S)_{-m-n} \neq 0\; \text{for some}\; i \le t\}\\
& = &\max\{m|\ \ext_S^i(M,S)_{-m-n} \neq 0\; \text{for some}\; i \ge n-t\}.\end{eqnarray*}
On the other hand, by [BCP, Proposition 1.2] (which is based on [E, Proposition 20.16]) we know that
$$\max\{m|\ \ext_S^i(M,S)_{-m-i} \neq 0\; \text{for some}\; i \ge n-t\} = \max\{b_i-i|\ i \ge n-t\}.$$
Hence (i) is immediate. \par
To prove (ii) we have to modify the proof of [E, Proposition 20.16] as follows. Put 
$$m' = \max\{b_i|\ i \ge n-t\}-n.$$ Let $i$ be any index $\ge n-t$. Then $F_i$ has no generators of degree $\ge m'+n+1$, so $F^*_i = \Hom_S(F_i,S)$ must be zero in degree $\le -m'-n-1$. Since $\ext_S^i(M,S)$ is the homology of the dual of the resolution of $M$ at $F^*_i$,  $\ext_S^i(M,S)_r = 0$ for $r \ge -m'-n-1$. Now let $i$ be the largest integer $\le n- t$ such that $b_i-n = m'$. Then $F^*_i$ has $S(m'+n)$ as a summand, whereas $F^*_{i+1}$ has no summand of the form $S(r)$ with $r \ge m'+n$. By the minimality of the resolution, the summand $S(m'+n)$ of $F^*_i$ must map to zero in  $F^*_{i+1}$. Moreover, nothing in $F^*_{i-1}$ can map on to the generator of $S(m'+n)$ in $F^*_i$, so it gives a nonzero class in $\ext_S^i(M,S)$ of degree $-m'-n$.  Thus,
$$\max\{m|\ \ext_S^i(M,S)_{-m-n} \neq 0\; \text{for some}\; i \ge n-t\} = m' = \max\{b_i|\ i \ge n-t\} - n,$$
which implies (ii).  \end{pf}  

Let $\syz_t(E)$ denote the $t$-th syzygy module of $M$ which is defined as the kernel of the map $F_t \to F_{t-1}$. There is the following relationships between the cohomological invariants  of $M$ and those of its syzygy modules.

\begin{Corollary} Let $M$ be an arbitrary graded $S$-module of finite type. Then\smallskip\par
{\rm (i)}  $\reg_t(M) = \reg(\syz_{n-t}(M))+n-t$,\smallskip\par
{\rm (ii)} $a^*_t(M) = a^*(\syz_{n-t}(M))$. \end{Corollary}

 \begin{pf}   Note that $0 \To F_r \To \ldots \To F_{n-t+1} \To F_{n-t} \To \syz_{n-t}(M)$ is a minimal free resolution of $\syz_{n-t}(M)$. Then applying Theorem 3.1 twice we get 
\begin{eqnarray*} \reg_t(M) & = & \max\{b_i-i|\ i \ge n-t\} =  \reg(\syz_{n-t}(M))+n-t,\\
a^*_t(M) & = & \max\{b_i|\ i \ge n-t\} - n = a^*(\syz_{n-t}(M)). \end{eqnarray*}  \end{pf} 

\begin{Corollary} Let $I$ be an arbitrary homogeneous ideal of $S$. Then \smallskip\par
{\rm (i)}  $\reg(\syz_t(I)) = \reg(\syz_t(\gin(I))$, \smallskip\par
{\rm (ii)} $a^*(\syz_t(I)) = a^*(\syz_t(\gin(I))$. \end{Corollary}

 \begin{pf}   This follows from Corollary 1.4 and Corollary 3.2.  \end{pf}  

Assume that $F_i = \oplus_jS(-j)^{\beta_{i,j}}.$ Then $\beta_{i,j}$ are called the Betti numbers of $M$. If $m = l$-$\reg(M) \ge (l+1)$-$\reg(M)$, then $\beta_{l,m+l}$ is called an extremal Betti number of $M$. This amounts to saying that $\beta_{l,m+l}  \neq 0$ and $\beta_{i,j+i} = 0$ for all $i \ge l$ and $j \ge m$. Hence, the extremal Betti numbers pinpoint  ``jumps" in the regularity of the successive syzygy modules. Bayer, Charalambous and Popescu [BCP, Theorem 1.6] proved that the extremal Betti numbers of a homogeneous ideal $I$ do not change when passing to a generic initial ideal of $I$. See also an alternate proof by Aramova and Herzog in [AH2] where they  extended this result to exterior algebras. Viewed in terms of local cohomology modules, an extremal Betti number of $M$ is the dimension of a graded piece of a local cohomology module of $M$ which corresponds to a ``jump" of the regularity $\reg_t(M)$.  

\section{Gr\"obner bases and reduction number} \smallskip

Let $I$ be an arbitrary homogeneous ideal of the polynomial ring $S = k[x_1,\ldots,x_n]$. Let $\frak m$ be the maximal graded ideal of the factor ring $S/I$. Let $J$ be an ideal of $S$ which contains $I$. Then ${\frak a} = J/I$ is a reduction of $\frak m$ if $S/J$ is of finite length, and $r_{\frak a}(S/I)$ is the largest non-vanishing degree of $S/J$. If $k$ is an infinite field, a reduction of $\frak m$ is minimal if and only if it is generated by $d$ elements, where $d = \dim S/I$. \par

Vasconcelos conjectured that $r(S/I) \le r(S/\In(I))$ [V, Conjecture 5.15]. Recently, Bresinsky and Hoa [BH, Theorem 12] proved this inequality for generic initial ideals. Inspired by their paper,  we will show that equality holds for the reverse lexicographic order. This will follow from the following observation.

\begin{Lemma}  Let $\In(I)$ denote the initial ideal of $I$ with respect to  the reverse lexicographic term order.  Assume that ${\frak a} = (I,x_n,\ldots,x_{n-d+1})/I$ is a minimal reduction of $\frak m$. Then ${\frak b} = (\In(I),x_n,\ldots,x_{n-d+1})/\In(I)$ is a minimal reduction of the maximal graded ideal of $S/\In(I)$ and $$r_{\frak a}(S/I) = r_{\frak b}(S/\In(I)).$$ \end{Lemma}
 
 \begin{pf}   Put $r_{\frak a}(S/I) = r$. Then $r$ is the largest non-vanishing degree of the factor ring $S/(I,x_n,\ldots,x_{n-d+1})$ and therefore of $S/\In(I,x_n,\ldots,x_{n-d+1})$ since these graded rings share the same Hilbert function. By [BS2, Lemma (2.2)] we have
$$\In(I,x_n,\ldots,x_{n-d+1}) = (\In(I),x_n,\ldots,x_{n-d+1}).$$
Hence $r$ is also the largest non-vanishing degree of  $S/(\In(I),x_n,\ldots,x_{n-d+1})$. Note that $\dim S/\In(I) = d$. Then we can conclude that $\frak b$ is a minimal reduction of the maximal graded ideal of $S/\In(I)$ and that $r_{\frak b}(S/\In(I)) = r$. \end{pf} 

The following result shows that generic minimal reductions always have the smallest reduction number.

\begin{Lemma}  Assume that $k$ is an infinite field. For a generic choice of linear forms $y_1,\ldots,y_d$,  ${\frak a} = (I,y_1,\ldots,y_d)/I$ is a minimal reduction of $\frak m$ with
$$r_{\frak a}(S/I) = r(S/I).$$ \end{Lemma}

 \begin{pf}   First observe that $\frak a$ is a minimal reduction of $\frak m$ if there is a number $r$ such that $(I,y_1,\ldots,y_d)_{r+1} = S_{r+1}$ and $r_{\frak a}(S/I)$ is the minimum of such numbers. It is clear that $(I,y_1,\ldots,y_d)_{r+1} = S_{r+1}$ if and only if $\dim_k(I,y_1,\ldots,y_d)_{r+1} = \dim_kS_{r+1}$. Let $y_i =\alpha_{i1}x_1 + \cdots + \alpha_{in}x_n,\; i = 1,\ldots,d,$ where $\alpha = (\alpha_{ij}) \in k^{dn}$. Then we can express the condition $\dim_k(I,y_1,\ldots,y_d)_{r+1} = \dim_kS_{r+1}$ as the non-vanishing of certain polynomial $f_r(u)$ at $\alpha$,  where $u = (u_{ij})$ is a family of $dn$ variables. Let $z_i = u_{i1}x_1 + \cdots + u_{in}x_n,\; i = 1,\ldots,d$. Put $S_u = k(u)[x_1,\ldots,x_n]$ and $I_u = IS_u$. If $r(S/I) = s$ and if $\frak a$ is any minimal reduction of $\frak m$ with  $r_{\frak a}(S/I) = s$ , then $f_s(\alpha) \neq 0$. Hence $f_s(u) \neq 0$. Therefore, ${\frak b} = (I_u,z_1,\ldots,z_d)/I_u$ is a minimal reduction of the maximal graded ideal of $S_u/I_u$ and $r_{\frak b}(S_u/I_u) \le s$. Put $r = r_{\frak b}(S_u/I_u)$. Then $f_{r-1}(u) = 0$ and $f_r(u) \neq 0$. Thus,  there is a non-empty open set of the space $k^{dn}$ such that if $\alpha \in U$ then $f_{r-1}(\alpha) = 0$ and $f_r(\alpha) \neq 0$. Hence $r_{\frak a}(S/I) = r$. So we obtain $r \ge r(S/I)$. Hence $r = r(S/I)$.  \end{pf}  

\begin{Theorem}  Let $\gin(I)$ denote the generic initial ideal of $I$ with repsect to  the reverse lexicographic order. Assume that $k$ is an infinite field. Then 
$$r(S/I) = r(S/\gin(I)).$$\end{Theorem}

 \begin{pf}   By Lemma 4.2 we may assume that the ideal ${\frak a} = (I,x_n,\ldots,x_{n-d+1})/I$, $d =\dim S/I$, is a minimal reduction of  $\frak m$ with $$r(S/I) = r_{\frak a}(S/I).$$
Let  ${\frak b} = (\gin(I),x_n,\ldots,x_{n-d+1})/\gin(I)$. By Lemma 4.1, $\frak b$ is a minimal reduction of the maximal graded ideal of $S/\gin(I)$ and
$$r_{\frak a}(S/I) =  r_{\frak b}(S/\gin(I)) \ge r(S/\gin(I)).$$
By [BH, Theorem 12] we know that $r(S/I) \le r(S/\gin(I))$, hence the conclusion.  \end{pf} 

\noindent{\em Remark. } The reduction number $r(S/\gin(I))$ can be easily computed. Bresinsky and Hoa [BH, Theorem 11] showed that $r(S/\gin(I))$ is the least number $r$ for which $x_{n-d}^{r+1} \in \gin(I)$. This  fact can be also deduced from Lemma 2.3. \smallskip

Now we will give an example showing that Theorem 4.3 does not hold for an arbitrary term order.\smallskip

\noindent{\em Example. } Let $S = k[x_1,x_2,x_3]$ and  $I = (x_1^2,x_1x_3-x_2^2)$. The ideal ${\frak a} = (I,x_3)/I$ is a minimal reduction of the maximal graded ideal of $S/I$ with
$$r_{\frak a}(S/I) = r(S/I) = 2.$$  It is not hard to check that $\gin(I) = (x_1^2,x_1x_2,x_1x_3^2,x_2^4)$ with respect to the lexicographic order. Let $\frak m$ denote the maximal graded ideal of $S/\gin(I)$. It is easy to verify that $H_{\frak m}^0(S/\gin(I)) = (x_1,x_2^4)/(x_1^2,x_1x_2^2,x_2^4,x_1x_3)$. Hence the largest non-vanishing degree of $H_{\frak m}^0(S/\gin(I))$ is 2 (that is the degree of $x_1x_3 \in \gin(I) $).  By [T1, Proposition 2.3 and Corollary 3.3] we have $\reg(S/\gin(I)) = \max\{2,r_{\frak b}(S/\gin(I))\}$
for any minimal reduction $\frak b$ of $\frak m$. If ${\frak b} = (\gin(I),x_3)/\gin(I)$, then $r_{\frak b}(S/\gin(I)) = 3$. Thus, $\reg(S/\gin(I)) = 3$. From this it follows that $r_{\frak b}(S/\gin(I)) = 3$ for any minimal reduction $\frak b$ of $\frak m$. Hence $$r(S/\gin(I)) = 3 > r(S/I).$$  

\section*{References}\smallskip

\noindent [AH1]  A.~Aramova and J.~Herzog, Koszul cycles and Eliahou-Kervaire type resolutions, J. Algebra 181 (1996), 347-370.\par

\noindent [AH2]  A.~Aramova and J.~Herzog, Almost regular sequences and Betti numbers, preprint.\par

\noindent [BCP] D.~Bayer, H.~Charalambous and S. Popescu, Extremal Betti numbers and applications to monomial ideals, preprint.\par

\noindent [BS1] D.~Bayer and M.~Stillman, A theorem on refining divison orders by the reverse lexicographic orders, Duke J. Math. 55 (1987), 321-328.\par

\noindent [BS2] D.~Bayer and M.~Stillman, A criterion for detecting $m$-regularity, Invent. Math. 87 (1987), 1-11.\par

\noindent [BH] H.~Bresinsky and L.T.~Hoa, On the reduction number of some graded algebras, Proc. Amer. Math. Soc. 127 (1999), 1257-1263.\par

\noindent [E] D.~Eisenbud, Commutative Algebra with a viewpoint toward Algebraic Geometry, Springer, 1994. \par

\noindent [EG] D.~Eisenbud and S.~Goto, Linear free resolutions and minimal multiplicities, J.~Algebra~88 (1984), 89-133.\par

\noindent [EK] S.~Eliahou and M.~Kervaire, Minimal resolutions of some monomial ideals, J. Algebra 129 (1990), 1-25.\par

\noindent [Ga] A. Galligo, Th\'eor\`eme de division et stabilit\'e en g\'eometrie analytique locale, Ann. Inst. Fourier 29 (1979), 107-184.\par

\noindent [GW] S. Goto and K. Watanabe, On graded ring I, J. Math. Soc. Japan 30 (1978), 179-213.\par

\noindent [O] A.~Ooishi, Castelnuovo regularity of graded rings and modules, Hiroshima Math. J. 12 (1982), 627-644. \par

\noindent [SCT] P.~Schenzel, N. V.~Trung, and N. T.~Cuong, Verallgemeinerte Cohen-Macaulay-Moduln, Math. Nachr. 85 (1978), 57-73.\par

\noindent [T1] N.~V.~Trung, Reduction exponent and degree bounds for the defining equations of a graded ring, Proc. Amer. Math. Soc. 102 (1987), 229-236.\par

\noindent [T2] N.~V.~Trung, The Castelnuovo regularity of the Rees algebra and the associated graded ring, Trans. Amer. Math. Soc. 350 (1998), 2813-2832. \par

\noindent [T3] N.~V.~Trung, The largest non-vanishing degree of graded local cohomology modules, J. Algebra 215 (1999), 481-499. \par

\noindent  [V] W.~Vasconcelos, Cohomological degrees of graded modules, in:  Six Lectures on Commutative Algebra, Progress in Mathematics 166, Birkh\"auser, Boston, 1998, 345-392.
\end{document}